\documentclass[review]{elsarticle}
\usepackage{graphicx,subfig}
\usepackage{lineno,hyperref}
\usepackage{amsmath,amssymb}
\usepackage{color}
\usepackage{float}
\modulolinenumbers[5]

\journal{Journal of Approximation Theory}

\begin{document}
\newcommand{\R}{\mathbb R}
\newcommand{\N}{\mathbb N}
\newcommand{\Z}{\mathbb Z}
\newcommand{\dd}{\mathrm{d}}
\newcommand{\lam}{\lambda}
\newcommand{\bfx}{\mathbf{x}}
\newcommand{\bfy}{\mathbf{y}}
\newcommand{\bfs}{\mathbf{s}}
\newcommand{\bfz}{\mathbf{z}}
\newcommand{\bfa}{\mathbf{a}}
\newcommand{\bfb}{\mathbf{b}}
\newcommand{\hK}{\widehat{K}}
\newcommand{\tK}{\widetilde{K}}
\newcommand{\cK}{\Check{K}}
\newcommand{\wA}{\widetilde{A}}
\newcommand{\hv}{\hat{v}}
\newcommand{\T}{\mathcal{T}}
\newcommand{\I}{\mathcal{I}}
\newcommand{\PP}{\mathcal{P}}
\newcommand{\Kab}{K_{ab}}
\newcommand{\Kabc}{K_{\alpha\beta\gamma}}
\newcommand{\Simp}{\mathbb S}
\newcommand{\bK}{\mathbf{K}}
\newcommand{\ux}{\underbar{x}}
\newcommand{\ox}{\overline{\mathrm{x}}}
\newcommand{\vvskip}{\vspace{5pt}}

\renewcommand{\figurename}{Figure~}
\renewcommand{\baselinestretch}{1.15}

\newtheorem{theorem}{Theorem}
\newtheorem{lemma}[theorem]{Lemma}
\newtheorem{definition}[theorem]{Definition}
\newtheorem{corollary}[theorem]{Corollary}

\begin{frontmatter}

\title{Error analysis of Lagrange interpolation on tetrahedrons}


\author[mymainaddress]{Kenta Kobayashi}
\ead{kenta.k@r.hit-u.ac.jp}

\author[mysecondaryaddress]{Takuya Tsuchiya\corref{mycorrespondingauthor}}
\cortext[mycorrespondingauthor]{Corresponding author}
\ead{tsuchiya@math.sci.ehime-u.ac.jp}

\address[mymainaddress]{Graduate School of Business Administration,
         Hitotsubashi University, Kunitachi, Japan}
\address[mysecondaryaddress]{Graduate School of Science and Engineering,
 Ehime University, Matsuyama, Japan}

\begin{abstract}
This paper describes the analysis of Lagrange interpolation errors on
tetrahedrons.  In many textbooks, the error analysis of Lagrange
interpolation is conducted under geometric assumptions such as shape regularity or
the (generalized) maximum angle condition.  In this paper,
we present a new estimation in which the error is
bounded in terms of the diameter and projected circumradius
of the tetrahedron.  Because we do not impose any geometric restrictions
on the tetrahedron itself, our error estimation may be applied to any
tetrahedralizations of domains including very thin tetrahedrons.
\end{abstract}

\begin{keyword}
Lagrange interpolation \sep tetrahedrons \sep
projected circumradius  \sep finite elements
\MSC[2010] 65D05 \sep  65N30
\end{keyword}

\end{frontmatter}

\linenumbers

\section{Introduction}\label{intro}
Lagrange interpolation on tetrahedrons and the associated error
estimates are important subjects in numerical analysis. In particular,
they are crucial in the error analysis of finite element methods.
Throughout this paper, $K \subset \R^3$ denotes a tetrahedron
with vertices $\bfx_i$, $i=1,\cdots,4$, and all tetrahedrons are assumed to
be closed sets. Let $\lam_i$ be the barycentric coordinates of a
tetrahedron with respect to $\bfx_i$. By definition,
we have $0 \le \lam_i \le 1$, $\sum_{i=1}^{4} \lam_i =1$.  Let $\N_0$
be the set of nonnegative integers and
$\gamma = (a_1,\cdots,a_{4}) \in \N_0^{4}$ be a multi-index.
Let $k$ be a positive integer. If $|\gamma| := \sum_{i=1}^{4}a_i = k$, then
${\gamma}/{k} := \left({a_1}/{k}, \cdots, {a_{4}}/{k}\right)$
can be regarded as a barycentric coordinate in $K$.
The set $\Sigma^k(K)$ of points on $K$ is defined by
\begin{equation*}
   \Sigma^k(K) := \left\{
   \frac{\gamma}{k} \in K \Bigm| |\gamma| = k, \;
     \gamma \in \N_0^{4} \right\}.
\end{equation*}
$\PP_k(K)$ is the set of all polynomials defined on $K$ whose
degree is at most $k$.  For a continuous function $v \in C^0(K)$, the
Lagrange interpolation $\I_K^k v \in \PP_k(K)$ of degree $k$ is defined as
\[
   v(\bfx) = (\I_K^k v)(\bfx), \quad \forall \bfx \in \Sigma^k(K).
\]
We attempt to obtain an upper bound of the error
$|v - \I_K^k v|_{m,p,K}$ for integers $0 \le m \le k$, where
$|\cdot|_{m,p,K}$ is the usual Sobolev semi-norm.
\footnote{
Note that, in this paper,  we deal with only Lagrange interpolation
defined with
``equally spaced points'' on simplices.  The accuracy of Lagrange
interpolation may be increased by moving points on simplices.
See \cite{Bos} and the references therein for this direction.}

Let us mention some well-known facts.
For now, let $K$ be a $d$-simplicial element ($d \ge 2$).  Let
$h_K := \mathrm{diam}K$ and $\rho_K$ be the diameter of its inscribed
sphere.  For the error analysis of Lagrange interpolation on simplicial
elements, many textbooks on finite element methods
\cite{Ciarlet, BrennerScott, ErnGuermond}
explain the following theorems \cite[Theorem~3.1.3, 3.1.4]{Ciarlet}.

Let $\hK$ be a reference element. Let $\varphi(\bfx) = A\bfx + \bfb$ be
an affine transformation that maps $\hK$ to $K$, where $A$ is
a $d \times d$ regular matrix and $\bfb \in \R^d$. 
Let $1 \le p \le \infty$.  A function
$v \in W^{2,p}(K)$ is pulled-back by $\varphi$ as
$\hat{v} := v \circ \varphi$, and error analysis is conducted on
$\hK$. 
\begin{theorem}
 We have $\|A\| \le h_K \rho_{\hK}^{-1}$, 
$\|A^{-1}\| \le h_{\hK} \rho_K^{-1}$, where $\|A\|$ denotes the matrix
norm of $A$ associated with the Euclidean norm of $\R^d$.
\end{theorem}

Let $k$, $m$ be integers such that $k \ge 1$
and $0 \le  m \le k$. 

\begin{theorem}[Shape-regularity]\label{thm:S-regularity}
Let $\sigma > 0$ be a constant. If $h_K/\rho_K \le \sigma$, then there
exists a constant $C = C(\hK,p,k,m)$ independent of $K$ such that,
for $v \in W^{k+1,p}(K)$,
\begin{align}
  |v - \I_K^k v|_{m,p,K} & \le C \|A\|^{k+1}\|A^{-1}\|^m
   |v|_{k+1,p,K} \notag \\
  & \le C \frac{h_K^{k+1}}{\rho_K^m} |v|_{k+1,p,K}
   \le (C\sigma^m) h_K^{k+1-m} |v|_{k+1,p,K}.
  \label{eq:standard-est}
\end{align}
\end{theorem}

For the case of triangles (that is, $d=2$), a triangle with vertices
$(0,0)^\top$, $(1,0)^\top$, and $(0,1)^\top$ is taken as the reference
triangle $\hK$.  Also, we may assume without loss
of generality that $K$ is the triangle with vertices $\bfx_1 = (0,0)^\top$,
$\bfx_2 = (\alpha,0)^\top$, $\bfx_3 = (\beta s, \beta t)^\top$, where
$\alpha \ge \beta > 0$, $s = \cos\theta$, $t = \sin\theta$, and
$0 < \theta < \pi$ is the inner angle of $K$ at $\bfx_1$.
Furthermore, we may also assume $\alpha \le |\bfx_2\bfx_3| = h_K$.
Then, we have $h_K/2 < \alpha \le h_K$ and $\pi/3 \le \theta < \pi$.
These assumptions imply that the affine transformation $\varphi$
can be written as $\varphi(\bfx) = A\bfx$ with the matrix
\begin{align}
  A = \begin{pmatrix}
      \alpha & \beta s \\
       0     & \beta t
      \end{pmatrix}.
   \label{mat-A}
\end{align} 
Set $t = \sin\theta = 1$, for example (that is, $K$ is a right triangle).
Then, $\|A\| = \alpha$ and $\|A^{-1}\| = 1/\beta$, and the inequalities
in \eqref{eq:standard-est} can be rearranged as
\begin{align}
  |v - \I_K^k v|_{m,p,K} \le C \frac{\alpha^{k+1}}{\beta^m}
   |v|_{k+1,p,K}  \le C
  \left(\frac{\alpha}{\beta}\right)^m h_K^{k+1-m} |v|_{k+1,p,K}.
   \label{ab-est}
\end{align}
Thus, one might consider that the ratio $\alpha/\beta$ should not be
too large, or $K$ should not be too ``flat.''  This consideration is
expressed as the \textit{minimum angle condition} \cite{Zlamal}, which
is equivalent to the shape-regularity condition for triangles.

It is known that shape-regularity is not an optimal condition on
the geometry of triangles.  If the maximum angle of a
triangle is less than a fixed constant $\theta_1 < \pi$, then the
estimation \eqref{eq:standard-est} holds after replacing
$C \sigma^m$ with a constant $C_2 = C_2(\theta_1,p,k,m)$.
This condition is known as the
\textit{maximum angle condition} \cite{BabuskaAziz}, \cite{Jamet}. 

For the case $d=2$, the present authors recently reported an error
estimation in terms of the circumradius of a triangle
\cite{KobayashiTsuchiya1}, \cite{KobayashiTsuchiya3},
\cite{KobayashiTsuchiya4}.
Let $R_K$ be the circumradius of a triangle $K$.
\begin{theorem}[Circumradius estimates]\label{thm:circumradius-est}
Let $K$ be an arbitrary triangle.  Then, for the $k$th-order Lagrange
interpolation $\I_K^k$ on $K$, the estimation
\begin{equation} \label{eq:circumradius-est}
   |v - \I_K^k v|_{m,p,K} \le C
   \left(\frac{R_K}{h_K}\right)^m h_K^{k+1-m} |v|_{k+1,p,K}
  =  C R_K^m h_K^{k+1-2m} |v|_{k+1,p,K}
\end{equation}
holds for any $v \in W^{k+1,p}(K)$,
where the constant $C=C(k,m,p)$ is independent of $K$.
\footnote{
Note that the circumradius estimation 
\eqref{eq:circumradius-est} with $m=k=1$ is closely related to the
definition of surface area \cite{KobayashiTsuchiya2}.}
\end{theorem}

In the proof of Theorem~\ref{thm:circumradius-est}, the main point
is that the matrix $A$ in \eqref{mat-A} can be decomposed as
\begin{align*}
   A = \wA D_{\alpha\beta}, \qquad
   \wA := \begin{pmatrix} 1 & s \\ 0 & t \end{pmatrix},
   \quad
   D_{\alpha\beta} := \begin{pmatrix} \alpha & 0 \\ 0 & \beta
       \end{pmatrix}.
\end{align*}
As pointed out by Babu\v{s}ka--Aziz \cite{BabuskaAziz} and the present
authors \cite{KobayashiTsuchiya4}, the linear transformation by
$D_{\alpha\beta}$ does not reduce the approximation property of Lagrange
interpolation at all, and only $\wA$ could make it ``bad.''
In the proof, the error of Lagrange interpolation
is bounded in terms of the singular values of
$\wA$ and $\wA^{-1}$ (or $\|\wA\|$ and $\|\wA^{-1}\|$).
Finally, these quantities (in terms of the singular values) are bounded
by means of the circumradius $R_K$
(see \cite{KobayashiTsuchiya3} for details).
Note that, setting $t = 1$ and $\beta = \alpha^2$ in
\eqref{eq:standard-est} (and \eqref{ab-est}), we realize that no matter
how much we try to analyze $\|A\|^{k+1}\|A^{-1}\|^m$, we cannot prove
Theorem~\ref{thm:circumradius-est}.   Although our novel strategy
is simple, it is more powerful than expected and has
been applied to several finite element error analysis
without the shape-regularity condition
\cite{KobayashiTsuchiya5, KobayashiTsuchiya6}.

For tetrahedrons, a generalized maximum angle condition was derived
by K\v{r}\'{i}\v{z}ek \cite{Krizek2} and Dur\'an \cite{Duran}.
Let $K$ be a tetrahedron, and let $\theta_1 < \pi$ be fixed.
If any inner angles of the faces and any dihedral angles between faces
of $K$ are less than or equal to $\theta_1$, then $K$ is said to
satisfy \textit{K\v{r}\'{i}\v{z}ek's maximum angle condition} with
$\theta_1$.  Then, the following theorem holds.
\begin{theorem}[Maximum angle condition]\label{thm:Maximum-angle}
Let $\pi/2 \le \theta_1 < \pi$ be a constant.
Suppose that a tetrahedron $K$ satisfies K\v{r}\'{i}\v{z}ek's maximum
angle condition with $\theta_1$.  Then, there exists a constant
$C = C(\theta_1,p)$ with $p > 2$ such that
\begin{equation}
  |v - \I_K^1 v|_{1,p,K} \le C  h_K |v|_{2,p,K},
   \qquad  \forall v \in W^{2,p}(K),
  \label{eq:maximum-est}
\end{equation}
where $C(\theta_1,p) = \mathcal{O}((p-2)^{-1/2})$ as $p \searrow 2$.
\end{theorem}

The aim of this paper is to extend Theorem~\ref{thm:Maximum-angle} and
derive a similar error estimation to that of Theorem~\ref{thm:circumradius-est}
for tetrahedrons under no specific geometric restrictions.

To extend the circumradius estimation \eqref{eq:circumradius-est}
to tetrahedrons, an immediate idea is to replace the circumradius
of a triangle with the radius of the circumsphere of a tetrahedron.
However, this idea can be immediately rejected by considering 
the tetrahedron $K$ with vertices
$\bfx_1 = (h,0,0)^\top$, $\bfx_2 = (-h,0,0)^\top$,
$\bfx_3 = (0,-h,h^\alpha)^\top$, and $\bfx_4 = (0,h,h^\alpha)^\top$
with $h > 0$ and $\alpha > 0$.
Setting $v_1(x,y,z) := x^2 - h^2 + h^{2-\alpha}z$,
we have that $\I_K^1 v_1 \equiv 0$, and a simple computation yields
$|v_1 - \I_K^1 v_1|_{1,\infty,K} = |v_1|_{1,\infty,K} \ge h^{2-\alpha}$
and $|v_1|_{2,\infty,K} = 2$.  Hence, if $\alpha > 2$,
an inequality such as \eqref{eq:circumradius-est} does
not hold for tetrahedrons, although the radius of the circumsphere of
the above $K$ converges to $0$ as $h \to 0$.  Tetrahedrons such as
$K$ above are called \textit{slivers} \cite{Sliver}.

Thus, we introduce the projected circumradius, also denoted by
$R_K$, of a tetrahedron $K$ in Section~\ref{section:Circumradius}.
We then obtain the following main theorem for tetrahedrons that is
fundamentally similar to \eqref{eq:circumradius-est}.

\begin{theorem}[Main Theorem]\label{thm:main}
Let $K$ be an arbitrary tetrahedron.  Let $h_K:=\mathrm{diam}K$ and
$R_K$ be the projected circumradius of $K$ defined by
\eqref{def:RK}.  Assume that $k$, $m$ are integers with $k \ge 1$,
$0 \le m \le k$, and take $p$ $(1 \le p \le \infty)$ as
\begin{gather}
  \begin{cases}
      2 < p \le \infty & \text{ if } k - m = 0, \\ 
      \frac{3}{2} < p \le \infty & \text{ if } k = 1, \; m = 0,\\
      1 \le p \le \infty & \text{ if } k \ge 2 \text{ and } \; k-m \ge 1.
   \end{cases}
   \label{p-cond}
\end{gather}
Then, for arbitrary $v \in W^{k+1,p}(K)$,
there exists a constant $C = C(k,m,p)$ independent of $K$ such that
\begin{equation}
\begin{split}
   |v - \I_K^k v|_{m,p,K} & \le C R_K^m h_K^{k+1-2m}|v|_{k+1,p,K} \\
 & = C \left(\frac{R_K}{h_K}\right)^m h_K^{k+1-m}|v|_{k+1,p,K}.
\end{split}
  \label{main-est}
\end{equation}
 \end{theorem}


Obviously, the error estimation \eqref{main-est} in terms of the
projected circumradius of tetrahedrons is completely different from most
prior error estimates obtained under the maximum angle condition. 
Additionally, as has been pointed out before, Theorem~\ref{thm:main}
cannot be proved by analyzing the inequality \eqref{eq:standard-est}.
As stated in Section~\ref{section:Circumradius}, there
exists a sequence $\{K_h\}$ of tetrahedrons with
$\lim_{h \to 0}R_{K_h} = 0$, whereas their maximum dihedral angles
approach $\pi$. Therefore, Theorem~\ref{thm:main} is an essential
extension of the prior error estimations, such as those in
Theorems~\ref{thm:S-regularity}, \ref{thm:Maximum-angle}.  Further
comments on these points are given in the concluding remarks (Section~6).

In Section~\ref{Sect:Prei}, we state some preliminary results used in
this paper.
In particular, we define the standard position and the projected
circumradius for tetrahedrons.  In Section~\ref{Sect:Squeezing}, we
recall definitions related to the quotient differences of functions
with multi-variables. We then reconfirm the Squeezing Theorem 
(Theorem~\ref{squeezingtheorem}) for tetrahedrons.
In Section~\ref{sect:singular}, we obtain the error
estimation of Lagrange interpolation in terms of the singular values of
a linear transformation.  In Section~\ref{sec:geometric}, we present
a geometric interpretation of the singular values of the linear
transformation, which completes the proof of the main theorem.
In the appendix, we give an interesting and  simple counter example
showing that the estimation \eqref{main-est} does not hold for the
case $p = 2$ and $k = m = 1$.

\section{Preliminaries}\label{Sect:Prei}
\subsection{Notation}
Let $d \ge 1$ be a positive integer and $\R^d$ be the $d$-dimensional
Euclidean space.  We denote the Euclidean norm of
$\bfx = (x_1,\cdots,x_d)^\top \in \R^d$ by $|\bfx|$. 
We always regard $\bfx \in \R^d$ as a column vector.
For a matrix $A$ and $\bfx \in \R^d$, $A^\top$ and $\bfx^\top$
denote their transpositions.

For $\delta = (\delta_1,...,\delta_d) \in \N_{0}^d$,
the multi-index $\partial^\delta$ of partial differentiation 
(in the sense of the distribution) is defined by
\[
    \partial^\delta = \partial_\bfx^\delta
    := \frac{\partial^{|\delta|}\ }
   {\partial x_1^{\delta_1}...\partial x_d^{\delta_d}}, \qquad
   |\delta| := \delta_1 + ... + \delta_d.
\]

Let $\Omega \subset \R^d$ be a (bounded) domain.  The usual Lebesgue
space is denoted by $L^p(\Omega)$ for $1 \le p \le \infty$.
For a positive integer $k$, the Sobolev space $W^{k,p}(\Omega)$ is
defined by
$\displaystyle
  W^{k,p}(\Omega) := 
  \left\{v \in L^p(\Omega) \, | \, \partial^\delta v \in L^p(\Omega), \,
   |\delta| \le k\right\}$.
For $1 \le p < \infty$, the norm and semi-norm of $W^{k,p}(\Omega)$ are defined
by
\begin{gather*}
  |v|_{k,p,\Omega} := 
  \biggl(\sum_{|\delta|=k} |\partial^\delta v|_{0,p,\Omega}^p
   \biggr)^{1/p}, \quad   \|v\|_{k,p,\Omega} := 
  \biggl(\sum_{0 \le m \le k} |v|_{m,p,\Omega}^p \biggr)^{1/p}
\end{gather*}
and $\displaystyle   |v|_{k,\infty,\Omega} := 
  \max_{|\delta|=k} \left\{\mathrm{ess}
   \sup_{\hspace{-5mm}\bfx \in\Omega}|\partial^\delta v(\bfx)|\right\}$,
 $\displaystyle   \|v\|_{k,\infty,\Omega} := 
  \max_{0 \le m  \le k} \left\{|v|_{m,\infty,\Omega}\right\}$.

\subsection{The imbedding theorem}
Let $1 < p \le \infty$. From Sobolev's imbedding theorem and Morry's
inequality, we have the continuous imbeddings
\begin{gather*}
   W^{2,p}(K) \subset C^{1,1-3/p}(K), \quad p > 3, \\
   W^{2,3}(K) \subset W^{1,q}(K) \subset C^{0,1-3/q}(K),
      \quad \forall q > 3, \\
  W^{2,p}(K) \subset W^{1,3p/(3-p)}(K) \subset C^{0,2-3/p}(K),
      \quad \frac{3}{2} < p < 3, \\
  W^{3,3/2}(K) \subset W^{2,3}(K) \subset W^{1,q}(K)
  \subset C^{0,1-3/q}(K), \quad \forall q > 3, \\
  W^{3,p}(K) \subset W^{2,3p/(3-p)}(K) \subset W^{1,3p/(3-2p)}(K)
  \subset C^{0,3-3/p}(K), \quad 1 < p < \frac{3}{2}.
\end{gather*}
Although Morry's inequality may not be applied, the continuous
imbedding $W^{3,1}(K)$ $\subset C^{0}(K)$ still holds.
For the imbedding theorem, see \cite{AdamsFournier},
\cite{Brezis}, and \cite{KJF}.  In the following, we assume that $p$
is such that the imbedding $W^{k+1,p}(K) \subset C^{0}(K)$
holds, that is,
\begin{align*}
     1 \le p \le \infty, \quad  \text{ if } k+1 \ge 3
      \quad \text{ and } \quad
    \frac{3}{2} < p \le \infty, \quad \text{ if } k+1 =2.
\end{align*}

\subsection{Standard position of tetrahedrons}
When considering tetrahedrons, it is convenient to define their ``standard
coordinates.''  Take any tetrahedron $K$ with vertices
$\bfx_i$, $i = 1, \dots, 4$.  The facet $B$ with vertices $\bfx_1$,
$\bfx_2$, $\bfx_3$ is regarded as the base of $K$. Let $\alpha$ and
$\beta$, $0 < \beta \le \alpha$, be the longest and shortest
lengths of the edges of $B$.   Without loss of generality, we assume
that $\bfx_1\bfx_2$ is the longest edge of $B$;
$|\bfx_1 - \bfx_2| = \alpha$.  Consider cutting $\R^3$
with the plane that contains the midpoint of the edge $\bfx_1\bfx_2$ and
is perpendicular to the vector $\bfx_1 - \bfx_2$.  Then, there exist
two cases: (i) $\bfx_3$ and $\bfx_4$ belong to the same half-space,
or (ii) $\bfx_3$ and $\bfx_4$ belong to different half-spaces.
If either of $\{\bfx_3, \bfx_4\}$ is on the partition
plane, we suppose that $\{\bfx_3, \bfx_4\}$ belong to the same
half-plane.  Let $\gamma := |\bfx_1 - \bfx_4|$.
Under appropriate rotation, translation, and reflection operations,
these situations can be written using the parameters
\begin{equation}
 \begin{cases}
  0 < \beta \le \alpha, \quad 0 < \gamma, \quad s_1^2 + t_1^2  = 1, \;
 s_1 > 0, \; t_1 > 0, \quad \beta s_1 \le \frac{\alpha}{2}, \\
 s_{21}^2 + s_{22}^2 + t_2^2 = 1, \; t_2 > 0, \quad 
  \gamma s_{21} \le \frac{\alpha}{2},
 \end{cases}
 \label{eq:tetra-param}
\end{equation}
as
\begin{subequations} \label{eq:position}
\begin{gather}
 \bfx_1 = (0,0,0)^\top, \; \bfx_2 = (\alpha,0,0)^\top, \;
  \bfx_4 = (\gamma s_{21}, \gamma s_{22}, \gamma t_2)^\top,\\
  \begin{cases}
   \bfx_3 = (\beta s_1, \beta t_1,0)^\top & \text{for the case (i)} \\
   \bfx_3 = (\alpha - \beta s_1, \beta t_1,0)^\top & \text{for the case (ii)}
  \end{cases}.
\end{gather}
\end{subequations}
Note that, in the above setting, we implicitly assume that
$\bfx_1$ and $\bfx_4$ belong to the same half-space.
We refer to the coordinates in \eqref{eq:position}
as the \textit{standard position} of a tetrahedron.  In the following,
we sometimes write $h_B := \alpha$.  Let $R_B$ be the circumradius
of $B$.

\subsection{The reference tetrahedrons}
\begin{figure}[tbhp]
\centering
\includegraphics[width=9cm]{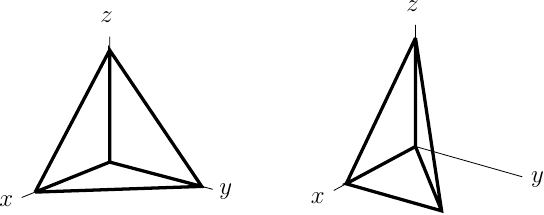}
\caption{The reference tetrahedrons $\hK$ and $\tK$.}
\label{ref_tetra}
\end{figure}

Because we have the two cases in the standard position of tetrahedrons,
we introduce two reference tetrahedrons to deal with the two cases.
Let $\hK$ and $\tK$ be tetrahedrons that have the 
following vertices (see Figure~\ref{ref_tetra}):
\begin{align*}
  \hK \text{ has the vertices } \; (0,0,0)^\top, \; (1,0,0)^\top, \;
  (0,1,0)^\top, \; (0,0,1)^\top, \\
  \tK \text{ has the vertices } \; (0,0,0)^\top, \; (1,0,0)^\top, \;
  (1,1,0)^\top, \; (0,0,1)^\top.
\end{align*}
These tetrahedrons are called the \textbf{reference tetrahedrons}.
In the proof of the main theorem, $\hK$ corresponds to the case (i) and
$\tK$ corresponds to the case (ii).
In this paper, we denote the reference tetrahedrons by 
$\bK$, that is, $\bK$ is either of $\{\hK,\tK\}$.

\subsection{The projected circurmradius of tetrahedrons}
\label{section:Circumradius}
\begin{figure}[tbhp]
\centering
\includegraphics[width=9cm]{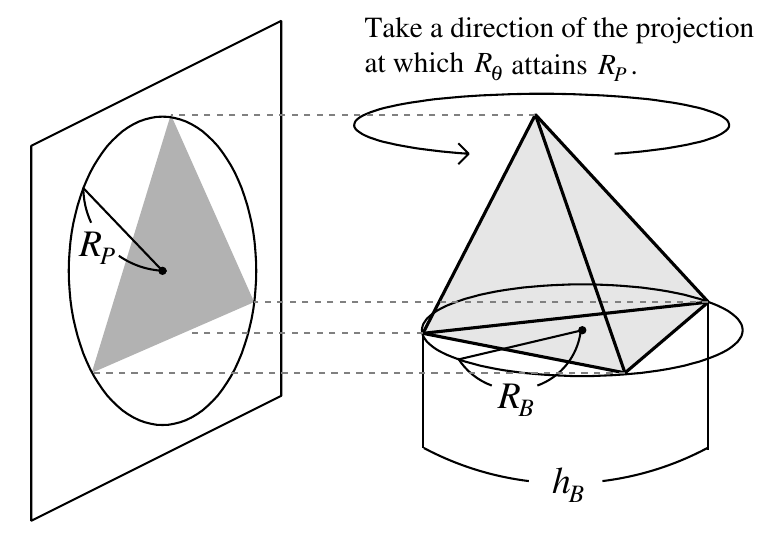}
\caption{The definition of the projected circumradius.}
\label{fig2}
\end{figure}

Suppose that a tetrahedron $K$ is at the standard position.
Let $\theta \in \R$ be such that
$-\frac{\pi}{2} \le \theta \le \frac{\pi}{2}$.
Let $\delta_\theta$ be the linear transformation defined by
the matrix
\[
   \begin{pmatrix}
      \cos\theta & -\sin\theta & 0 \\
          0 & 0 & 0 \\
          0 & 0 & 1
   \end{pmatrix}
   = \begin{pmatrix}
      1 & 0 & 0 \\
      0 & 0 & 0 \\
      0 & 0 & 1
   \end{pmatrix}
   \begin{pmatrix}
      \cos\theta & -\sin\theta & 0 \\
      \sin\theta & \cos\theta & 0 \\
          0 & 0 & 1
   \end{pmatrix}.
\]
That is, $\delta_\theta$ is a composite transformation of rotation
about the $z$-axis with angle $\theta$ and projection to the $xz$-plane.
The image $\delta_\theta(K)$ is a triangle on the $xz$-plane.  Let
$R_\theta$ be the circumradius of $\delta_\theta(K)$.  Define
\[
   R_P := \max_{\theta \in [-\pi/2,\pi/2]} R_\theta.
\]
The \textit{projected circumradius} $R_K$ of a tetrahedron $K$ is
defined by
\begin{equation}\label{def:RK}
   R_K := \min_{B} \frac{R_BR_P}{h_B},
\end{equation}
where the minimum is taken over all the facets of $K$.
See Figure~\ref{fig2} for a graphical image of the projected
circumradius.

Let $\delta$, $\eta >0$ be positive numbers with
$1 < \delta < \eta < 1 + \delta$.
Let $K$ be the tetrahedron with vertices
$\bfx_1 := (0,0,0)^\top$, $\bfx_2 := (h,0,0)^\top$ ,
$\bfx_3 := (h^\delta,h^\eta, 0)$, and $\bfx_4 := (0,0,h)^\top$.
Then, taking the triangle with vertices $\bfx_1$, $\bfx_2$, $\bfx_3$
as $B$, we realize that $R_B = \mathcal{O}(h^{1+\delta-\eta})$,
$h_B = h$, $R_P = \mathcal{O}(h)$, and
$R_K = \mathcal{O}(h^{1+\delta-\eta})$.
Thus, $R_K \to 0$ as $h \to 0$ whereas the maximum angle of $K$
approaches to $\pi$.  This example shows that
the condition ``$R_K \to 0$'' is more general than the maximum angle
condition.

\subsection{The setting of error estimation}
We define the set $\T_p^{k}(K) \subset W^{k+1,p}(K)$ by
\begin{gather*}
   \T_p^{k}(K) := \left\{v \in W^{k+1,p}(K) \Bigm|
     v(\bfx) = 0, \; \forall \bfx \in \Sigma^k(K) \right\}.
\end{gather*}
For Lagrange interpolation $\I_K^k(v)$, it is clear from the definition that
\[
     v - \I_K^k v \in \T_p^{k}(K), \quad \forall 
     v \in W^{k+1,p}(K).
\]
For an integer $m$ such that $0 \le m \le k$, $B_p^{m,k}(K)$
is defined by
\begin{gather*}
    B_p^{m,k}(K) := \sup_{v \in \T_p^{k}(K)}
    \frac{|v|_{m,p,K}}{|v|_{k+1,p,K}}.
\end{gather*}
Note that we have
\begin{gather*}
 B_p^{m,k}(K) = \inf\left\{C ; 
   |v - \I_K^k v|_{m,p,K} \le C |v|_{k+1,p,K}, \
   \forall v \in W^{k+1,p}(K)\right\},
\end{gather*}
that is, $B_p^{m,k}(K)$ is the \textit{best} constant $C$ for the
error estimation
\[
   |v - \I_K^k v|_{m,p,K} \le C |v|_{k+1,p,K}, \qquad 
    \forall v \in W^{k+1,p}(K).
\]
Therefore, we try to obtain an upper bound of $B_p^{m,k}(K)$ in terms of
geometric quantities of $K$.

\section{The Squeezing Theorem}\label{Sect:Squeezing}
Let $a$, $b \in\R$ be such that $0 < a \le 1$ and $0 < b$.
We then define the \textit{squeezing map} $sq_1^{ab}:R^3 \to \R^3$ by
\[
  sq_1^{ab}(x,y,z) := (x, a y, b z)^\top,
   \quad (x,y,z)^\top \in \R^3.
\]
Let $\Kab := sq_1^{ab}(\bK)$.  In this section, we will prove the
following theorem that is called the \textit{Squeezing Theorem}.

\begin{theorem}\label{squeezingtheorem}
Let $k$ and $m$ be integers with $k \ge 1$ and $0 \le m \le k$.
Let $p$ be taken as \eqref{p-cond}. Then, there exists a constant
$C_{k,m,p}$ depending only on $k$, $m$, $p$, but independent of
$a$ and $b$ such that
\begin{gather*}
   B_p^{m,k}(\Kab) := \hspace{-0.2cm}
   \sup_{v \in \T_p^{k}(\Kab)}
  \frac{|v|_{m,p,\Kab}}{|v|_{k+1,p,\Kab}}
   \le \max\{1,b^{k+1-m}\}C_{k,m,p}.
\end{gather*}
\end{theorem}
Although the proof of Theorem~\ref{squeezingtheorem} is very similar to
that of \cite[Theorem~1.3]{KobayashiTsuchiya4}, we provide a sketch here
for readers' convenience.  Note that, in \eqref{p-cond}, the restriction 
$2 < p \le \infty$ for the case $k=m$ comes from the continuity of
the trace operator $\gamma:W^{1,p}(\bK) \ni v \mapsto v|_S \in L^1(S)$,
where $S \subset \bK$ is a non-degenerate segment (see
\cite[Section~3]{KobayashiTsuchiya4}).
Using the counter-example given by Shenk \cite{Shenk}, we find that
this restriction cannot be improved.
\footnote{See the simpler counter-example given in Appendix.}

\subsection{Difference quotient of functions with three variables}
To prove Theorem~\ref{squeezingtheorem}, we briefly recall the
definitions of difference quotients of multi-variable functions, the
rectangular parallelepiped $\square_\gamma^\delta$ defined by the
lattice points $\bfx_\gamma$, $\Delta^\delta\bfx_\gamma$ in $\bK$, and
the integral $\int_{\square_\gamma^\delta} v$ introduced in
\cite[Section~2]{KobayashiTsuchiya4}. Note that the integral of
difference quotients have already been introduced and used for the
finite element error analysis on anisotropic meshes
\cite{ApelDobrowolski, ApelLube}.
For the definition of difference
quotients of functions and their basic properties, 
readers are referred to the textbooks
\cite{Atkinson, Krylov, Yamamoto1}.

For a positive integer $k$, $X^k$ is the set of lattice points
defined by
\begin{align*}
   X^k & := \left\{ \bfx_{\gamma} := \frac{\gamma}{k}
     \in \R^3 \biggm|  \gamma \in \N_0^3 \right\},
\end{align*}
where $\gamma/k=(a_1/k,a_2/k,a_d/k)$ is understood as the
coordinate of a point in $\R^3$.
For $\bfx_{\gamma} \in X^k$ and a multi-index $\delta \in \N_0^3$,
we define the correspondence $\Delta^\delta$
between nodes by
$\Delta^\delta \bfx_{\gamma} := \bfx_{\gamma+\delta} = (\gamma+\delta)/k$.

For two multi-indexes  $\eta=(m_1,m_2,m_3)$,
$\delta=(n_1,n_2,n_3)$,  $\eta \le \delta$ means that
$m_i \le n_i$ $(i=1,2,3)$.
Also, $\delta\cdot\eta$ and $\delta !$ are defined by
$\delta\cdot\eta := \sum_{i=1}^3 m_i n_i$ and
 $\delta ! := n_1! n_2! n_3!$, respectively.  Suppose that, for
$\gamma, \delta \in \N_0^3$, both $\bfx_\gamma$ and
$\Delta^\delta \bfx_\gamma$ belong to $\bK$.  Then, we define the
\textit{difference quotients} for $f \in C^0(\bK)$ by
\begin{align*}
 f^{|\delta|}[\bfx_{\gamma},\Delta^\delta \bfx_{\gamma}] :=
  k^{|\delta|}\sum_{\eta \le \delta} 
\frac{(-1)^{|\delta|-|\eta|}}{\eta!(\delta-\eta)!}
    f(\Delta^\eta \bfx_{\gamma}).
\end{align*}

If $f \in C^k(\bK)$, the difference quotient
$f^{|\delta|}[\bfx_{\gamma},\Delta^\delta \bfx_{\gamma}]$ is written as
an integral of $f$.  For example, we have
{\allowdisplaybreaks
\begin{align*}
   f^{1}[\bfx_{(l,q,r)},\Delta^{(0,1,0)} \bfx_{(l,q,r)}]
   & = k(f(\bfx_{(l,q+1,r)}) - f(\bfx_{(l,q,r)})) \\
   & = \int_{0}^{1} \partial^{(0,1,0)}
    f\left(\frac{l}{k},\frac{q}{k} + \frac{w_1}{k},
        \frac{r}{k} \right) \dd w_1,
\end{align*}
\begin{align*}
  f^{1}[\bfx_{(l,q,r)}, \Delta^{(0,2,0)} \bfx_{(l,q,r)}]
 & = \frac{k^2}{2}(f(\bfx_{(l,q+2,r)}) - 2 f(\bfx_{(l,q+1,r)}) + f(\bfx_{(l,q,r)})) \\
 & = \int_{0}^{1}\int_0^{w_1} \partial^{(0,2,0)}
  f\left(\frac{l}{k},\frac{q}{k} + \frac{1}{k}(w_1 + w_2), \frac{r}{k}
   \right) \dd w_2 \dd w_1.
\end{align*}
\begin{gather*}
 \hspace{-8truecm}
f^{s}[\bfx_{(l,p,r)},\Delta^{(0,s,0)} \bfx_{(l,q,r)}]  \\
\hspace{1truecm}
= \int_{0}^{1}\int_0^{w_1}\cdots \int_0^{w_{s-1}} \partial^{(0,s,0)}
  f\left(\frac{l}{k},\frac{q}{k} + \frac{1}{k}(w_1 + \cdots + w_s),  \frac{r}{k}
   \right) \dd w_s \cdots \dd w_2 \dd w_1.
\end{gather*}
}

To provide a concise expression for the above integral, we introduce
the $s$-simplex
\begin{gather*}
   \Simp_s := \left\{(x_1,\cdots,x_s)^\top \in \R^s \mid 
   x_i \ge 0, \ 0 \le x_1 + \cdots + x_s \le 1 \right\},
\end{gather*}
and the integral of $g \in L^1(\Simp_s)$ on $\Simp_s$ is defined by
\begin{gather*}
  \int_{\Simp_s} g(w_1,\cdots,w_k) \dd\mathbf{W_s}
  := \int_{0}^{1}\int_0^{w_1}\cdots \int_0^{w_{s-1}} 
  g(w_1, \cdots, w_s)  \dd w_s \cdots \dd w_2 \dd w_1,
\end{gather*}
where $\dd \mathbf{W_s} := \dd w_1 \cdots \dd w_s$.
Then, $f^{s}[\bfx_{(l,q)},\Delta^{(0,s,0)} \bfx_{(l,q)}]$ becomes
\begin{align*}
f^{s}[\bfx_{(l,q,r)},\Delta^{(0,s,0)} \bfx_{(l,q,r)}]
  & = \int_{\Simp_s} \partial^{(0,s,0)}
  f\left(\frac{l}{k},\frac{q}{k} + \frac{1}{k}(w_1 + \cdots + w_s),
    \frac{r}{k} \right) \dd \mathbf{W_s}.
\end{align*}
For a general multi-index $(t,s,m)$, we can write
\begin{gather*}
 f^{t+s+m} [\bfx_{(l,q,r)}, \Delta^{(t,s,m)} \bfx_{(l,q,r)}] 
  = \int_{\Simp_s}\int_{\Simp_t} \int_{\Simp_m}\partial^{(t,s,m)}
  f\left(\mathbf{Z_t}, \mathbf{W_s}, \mathbf{Y_m}\right)
    \dd \mathbf{Z_t} \dd \mathbf{W_s} \dd \mathbf{Y_m}, \\
  \mathbf{Z_t} := \frac{l}{k} + \frac{1}{k}(z_1 + \cdots + z_t), \;
   \dd \mathbf{Z_t} := \dd z_1 \cdots \dd z_t, \quad
  \mathbf{W_s} := \frac{q}{k} + \frac{1}{k}(w_1 + \cdots + w_s),
     \\
  \mathbf{Y_m} := \frac{r}{k} + \frac{1}{k}(y_1 + \cdots + y_m), \quad
    \dd \mathbf{Y_m} := \dd y_1 \cdots \dd y_m.
\end{gather*}

Let $\square_{\gamma}^\delta$ be the rectangular parallelepiped defined by 
$\bfx_{\gamma}$ and $\Delta^\delta \bfx_{\gamma}$ as the diagonal
points. If $\delta=(t,s,0)$ or $(0,s,0)$,  $\square_{\gamma}^\delta$
degenerates to a rectangle or a segment.  For $v \in L^1(\hK)$
and $\square_{\gamma}^\delta$ with $\gamma=(l,q,r)$, we denote the integral as
\begin{equation*}
   \int_{\square_{\gamma}^{(t,s,m)}} v :=
   \int_{\Simp_s}\int_{\Simp_t} \int_{\Simp_m}
  v\left(\mathbf{Z_t}, \mathbf{W_s}, \mathbf{Y_m}
   \right) \dd \mathbf{Z_t} \dd \mathbf{W_s} \dd \mathbf{Y_m}.
\end{equation*}
If $\square_{\gamma}^\delta$ degenerates to a rectangle or a segment, the integral is
understood as an integral on the rectangle or on the segment.  By this notation, the
difference quotient $f^{|\delta|}[\bfx_{\gamma},\Delta^{\delta} \bfx_{\gamma}]$ is
written as 
\begin{align*}
 f^{|\delta|}[\bfx_{\gamma},\Delta^{\delta} \bfx_{\gamma}] 
  = \int_{\square_{\gamma}^{\delta}} \partial^{\delta} f.
\end{align*}
Therefore, if $u \in \T_p^k(\hK)$, then we have
\begin{align}
 0 = u^{|\delta|}[\bfx_{\gamma},\Delta^{\delta} \bfx_{\gamma}] 
  = \int_{\square_{\gamma}^{\delta}} \partial^{\delta} u, \qquad
  \forall \square_{\gamma}^{\delta} \subset \bK.
  \label{tomoko}
\end{align}

\subsection{A proof of Theorem~\ref{squeezingtheorem}}
The set
$\Xi_p^{\delta,k} \subset W^{k+1-|\delta|,p}(\bK)$ is then defined by
\begin{align*}
  \Xi_p^{\delta,k} & := \left\{ v \in W^{k+1-|\delta|,p}(\bK) \Bigm| 
    \int_{\square_{\gamma}^{\delta}} v = 0,\quad \forall
   \square_{\gamma}^{\delta} \subset \bK \right\}.
\end{align*}
Note that $u \in \T_p^{k}(\bK)$ implies
$\partial^\delta u \in \Xi_p^{\delta,k}$ by \eqref{tomoko}.

\begin{lemma}\label{lem32} 
We have $\Xi_p^{\delta,k} \cap \mathcal{P}_{k-|\delta|} = \{0\}$.
That is, if $q \in \mathcal{P}_{k-|\delta|}$ belongs to $\Xi_p^{\delta,k}$,
then $q=0$.
\end{lemma}

\begin{figure}[htbp]
\centering
\subfloat{\includegraphics[width=4.05cm]{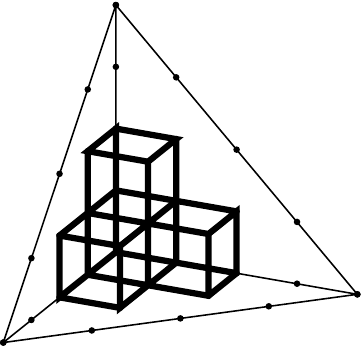}} \qquad
\subfloat{\includegraphics[width=2.7cm]{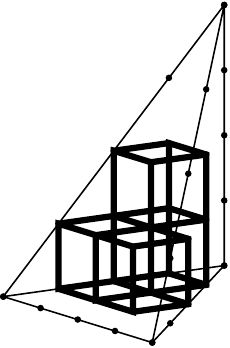}} \qquad
\subfloat{\includegraphics[width=2.7cm]{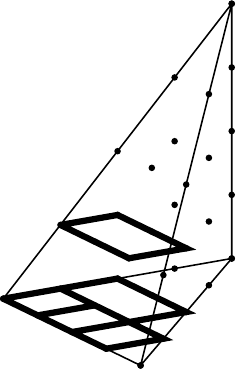}}
\caption{The four cubes and four rectangles in $\bK$.}
\label{cubes_tetra}
\end{figure}

\textit{Proof.}
 Note that
$\mathrm{dim} \mathcal{P}_{k-|\delta|} = \#\{\square_{\gamma}^\delta \subset \bK \}$.
For example, if $k=4$ and $|\delta|=3$, then
$\mathrm{dim}\mathcal{P}_1 = 4$.  This corresponds to the fact that,
in $\bK$, there are four cubes of size $1/4$  for $\delta=(1,1,1)$ and
there are four rectangles of size $1/2 \times 1/4$ for $\delta=(1,2,0)$.
All their vertices (corners) belong to $\Sigma^4(\bK)$
(see Figure~\ref{cubes_tetra}). Now, suppose that
$v \in  \mathcal{P}_{k-|\delta|}$ satisfies
$\int_{\square_{\gamma}^{\delta}} q = 0$ for all
$\square_{\gamma}^{\delta} \subset \bK$.  These conditions are linearly
independent and determine $q = 0$ uniquely. \hfill $\square$

\vspace{0.3cm}
The constant $A_p^{\delta,k}$ is defined by
\begin{align*}
   A_p^{\delta,k} := \sup_{v \in \Xi_p^{\delta,k}} \frac{|v|_{0,p,\bK}}
  {|v|_{k+1-|\delta|,p,\bK}}.
\end{align*}
The following lemma is an extension of \cite[Lemma~2.1]{BabuskaAziz}.

\begin{lemma}\label{lem33}
Let $p$ be such that $2 < p \le \infty$ if $k+1-|\delta|=1$ or
$1 \le p \le \infty$ if $k+1-|\delta|\ge 2$.  We then have
$A_p^{\delta,k} < \infty$.
\end{lemma}

\textit{Proof.}
See the proof of \cite[Lemma~3.3]{KobayashiTsuchiya4}.
\hfill $\square$

\vspace{0.3cm}
\textit{Proof of} Theorem~\ref{squeezingtheorem}.
First, let $1 \le p < \infty$ and $1 \le m \le k$.
For a multi-index $\gamma = (n_1, n_2, n_3) \in \N_0^3$ and a real $t \neq 0$,
set $(a, b)^{\gamma t} := a^{n_2 t}b^{n_3 t}$.  Take an arbitrary 
$v \in \T_p^k(\Kab)$ and pull it back to
$u := v\circ sq_1^{ab} \in \T_p^k(\bK)$.   We then have
{\allowdisplaybreaks
\begin{align*}
  \frac{|v|_{m,p,\Kab}^p}{|v|_{k+1,p,\Kab}^p}
   & = \frac{\sum_{|\gamma| = m} \frac{m!}{\gamma!}
      (a,b)^{-\gamma p}\left|\partial^{\gamma}u\right|_{0,p,\bK}^p}
      {\sum_{|\delta|= k+1} \frac{(k+1)!}{\delta!}
     (a,b)^{- \delta p} \left|\partial^{\delta} u \right|_{0,p,\bK}^p} \\
  & = \frac{\sum_{|\gamma| = m}\frac{m!}{\gamma!}(a,b)^{-\gamma p}
         \left|\partial^{\gamma}u \right|_{0,p,\bK}^p }
       {\sum_{|\gamma| = m} \frac{m!}{\gamma!}(a,b)^{-\gamma p}
       \left(\sum_{|\eta| = k+1-m} \frac{(k+1-m)!}
             {\eta!(a,b)^{\eta p} }\left|\partial^{\eta}
                 (\partial^{\gamma}u)\right|_{0,p,\bK}^p\right)} \\
  & \le \frac{\max\{1,b^{(k+1-m)p}\}
   \sum_{|\gamma| = m}\frac{m!}{\gamma!}(a,b)^{-\gamma p}
         \left|\partial^{\gamma}u \right|_{0,p,\bK}^p }
       {\sum_{|\gamma| = m} \frac{m!}{\gamma!}(a,b)^{-\gamma p}
       \left(\sum_{|\eta| = k+1-m} \frac{(k+1-m)!}
             {\eta! }\left|\partial^{\eta}
                 (\partial^{\gamma}u)\right|_{0,p,\bK}^p\right)} \\
   & = \max\{1,b^{(k+1-m)p}\}
\frac{\sum_{|\gamma| = m} \frac{m!}{\gamma!}(a,b)^{-\gamma p}
         \left|\partial^{\gamma}u \right|_{0,p,\bK}^p}
       {\sum_{|\gamma| = m} \frac{m!}{\gamma!}
      (a,b)^{-\gamma p}
        \left|\partial^{\gamma}u \right|_{k+1-m,p,\bK}^p} \\
    & \le \max\{1,b^{(k+1-m)p}\}\max_{|\gamma|=m} 
    \left(A_p^{\gamma,k}\right)^p
   = \left(\max\{1,b^{k+1-m}\}C_{k,m,p}\right)^p,
\end{align*}
}
where $C_{k,m,p} := \max_{|\gamma|=m} A_p^{\gamma,k}$.
Proofs for the cases $1 \le p < \infty$ with $m=0$ and
$p = \infty$ with $0 \le m \le k$ are very similar.  See
the proof of \cite[Theorem~1.3]{KobayashiTsuchiya4}.  
\hfill $\square$

\vspace{0.3cm}
We now generalize the squeezing map.
Let $\alpha$, $\beta$, and $\gamma \in\R$ be such that
$0 < \beta \le \alpha$ and $0 < \gamma$.
We then define the \textit{squeezing map}
$sq_2^{\alpha\beta\gamma}:\R^3 \to \R^3$ by
\[
  sq_2^{\alpha\beta\gamma}(x,y,z) := (\alpha x, \beta y, \gamma z)^\top,
   \quad (x,y,z)^\top \in \R^3.
\]
Let $\Kabc := sq_2^{\alpha\beta\gamma}(\bK)$.  Note that
$G_{\alpha} := sq_2^{\alpha\alpha\alpha}$ is a similar transformation.
Let $\Omega \subset \R^3$ be a domain and an arbitrary function 
$v \in W^{k,p}(\Omega)$ be pulled-back to
$u := v \circ G_\alpha \in W^{k,p}(\Omega_\alpha)$ with
$\Omega_\alpha := G_{1/\alpha}(\Omega)$.
It is straightforward to check that 
\[
  |v|_{k,p,\Omega} = \alpha^{3/p-k}|u|_{k,p,\Omega_\alpha}.
\]

Now, take an arbitrary function $v \in \T_p^k(\Kabc)$ and define
$u := v\circ sq_2^{\alpha\beta\gamma}$.  Then,
$u = u_1 \circ sq_1^{\frac{\beta}{\alpha},\frac{\gamma}{\alpha}}$ with
$u_1 := v \circ G_\alpha$, because
$sq_2^{\alpha\beta\gamma} = G_\alpha\circ sq_1^{\frac{\beta}{\alpha},\frac{\gamma}{\alpha}}$.
Therefore, it follows from Theorem~\ref{squeezingtheorem} that
\begin{align*}
   \frac{|v|_{m,p,\Kabc}}{|v|_{k+1,p,\Kabc}} = \alpha^{k+1-m}
    \frac{|u_1|_{m,p,K_{\frac{\beta}{\alpha},\frac{\gamma}{\alpha}}}}
   {|u_1|_{k+1,p,K_{\frac{\beta}{\alpha},\frac{\gamma}{\alpha}}}}
  & \le \alpha^{k+1-m} \max \left\{1,
   \left(\frac{\gamma}{\alpha}\right)^{k+1-m}\right\} C_{k,m,p} \\
  & = \left(\max\{\alpha,\gamma\}\right)^{k+1-m}C_{k,m,p}.
\end{align*}
Hence, we have derived the following corollary.

\begin{corollary}\label{cor:squeezing}
 Let $\alpha$, $\beta$, and $\gamma \in\R$ be such that
$0 < \beta \le \alpha$ and $0 < \gamma$.  
Let $\Kabc := sq_2^{\alpha\beta\gamma}(\bK)$. 
Assume that $k \ge 1$, $0 \le m \le k$, and $p$ is taken as
\eqref{p-cond}.   We then have
\begin{align*}
  B_p^{m,k}(\Kabc) := \hspace{-0.2cm} \sup_{v \in \T_p^{k}(\Kabc)}
  \frac{|v|_{m,p,\Kabc}}{|v|_{k+1,p,\Kabc}}
  \le \left(\max\{\alpha,\gamma\}\right)^{k+1-m}C_{k,m,p}. \label{extension-3d}
\end{align*}
\end{corollary}

\section{Error estimates of Lagrange interpolation on general
 tetrahedrons}
\label{sect:singular}
In this section, we obtain an error estimation for Lagrange interpolation
on general tetrahedrons.  To this end, we apply the method developed in
\cite{KobayashiTsuchiya3}.

Recall that an arbitrary tetrahedron $K$ is written as
\eqref{eq:position} with parameters \eqref{eq:tetra-param}.  First, we
confirm that $K$ is obtained from the reference tetrahedron $\bK$ by an
affine linear transformation.  Define the matrices $\widehat{A}$,
$\widetilde{A}$, $G \in GL(3,\R)$ by
\begin{align*}
   \widehat{A} := \begin{pmatrix}
        1  & s_{1} & s_{21} \\
        0  & t_{1} & s_{22} \\
        0  & 0     & t_2    
       \end{pmatrix}, \quad
   \widetilde{A} := \begin{pmatrix}
        1  & - s_{1} & s_{21} \\
        0  & t_{1} & s_{22} \\
        0  & 0     & t_2    
       \end{pmatrix}, \quad
  G := \begin{pmatrix}
        \alpha  & 0 & 0 \\
        0  & \beta & 0 \\
        0  & 0     & \gamma    
       \end{pmatrix}.
\end{align*}
We then have $K = \widehat{A}G(\hK)$ for case (i) or
$K = \widetilde{A}G(\tK)$ for case (ii), that is,
$K = \widehat{A}(\Kabc)$ or $K = \widetilde{A}(\Kabc)$.
Note that $\widehat{A}$ and $\widetilde{A}$ can be decomposed as
$\widehat{A} = X\widehat{Y}$ and $\widetilde{A} = X\widetilde{Y}$
with
\begin{align*}
  X := \begin{pmatrix}
        1 & 0 & s_{21} \\
	0 & 1 & s_{22} \\
        0 & 0 & t_2
       \end{pmatrix}, \qquad
  \widehat{Y} := \begin{pmatrix}
        1 & s_1 & 0 \\
	0 & t_1 & 0 \\
        0 & 0 & 1
       \end{pmatrix}, \qquad
  \widetilde{Y} := \begin{pmatrix}
        1 & -s_1 & 0 \\
	0 & t_1 & 0 \\
        0 & 0 & 1
       \end{pmatrix},
\end{align*}
respectively.  We consider the singular values of $\widehat{A}$,
$\widetilde{A}$, $X$, $\widehat{Y}$, and $\widetilde{Y}$.
A straightforward computation yields 
\begin{gather*}
  \det\left(X^\top X - \mu I\right) 
    = (1 - \mu)\left(\mu^2 - 2\mu + t_2^2\right), \\
  \det\left(\widehat{Y}^\top \widehat{Y} - \mu I\right) 
  = \det\left(\widetilde{Y}^\top \widetilde{Y} - \mu I\right) 
    = (1 - \mu)\left(\mu^2 - 2\mu + t_1^2\right),
\end{gather*}
the eigenvalues of which are
$\mu = 1$, $1 \pm \sqrt{1 - t_i^2} = 1 \pm \bfs_i$, $i = 1, 2$, where
$\bfs_1 := |s_1|$ and $\bfs_2 := \sqrt{s_{21}^2 + s_{22}^2}$.
Therefore, for $\bfa \in \R^3$, we have
\begin{gather*}
  (1 - \bfs_2)|\bfa|^2 \le |X\bfa|^2 \le (1 + \bfs_2)|\bfa|^2, \\
  (1 - \bfs_1)|\bfa|^2 \le |Z\bfa|^2 \le (1 + \bfs_1)|\bfa|^2,
  \qquad Z = \widehat{Y} \text{ or } Z = \widetilde{Y},  \\
 \prod_{i=1}^2(1 - \bfs_i)|\bfa|^2 \le |V\bfa|^2 \le 
 \prod_{i=1}^2(1 + \bfs_i)|\bfa|^2, \qquad
 V = \widehat{A} \text{ or } V = \widetilde{A}.
\end{gather*}

Let $K := V(\Kabc)$, where $V = \widehat{A}$ or $V = \widetilde{A}$.
A function $v\in W^{r,p}(K)$ is pulled-back
to a function $u \in W^{r,p}(\Kabc)$ by $u(\bfx) := v(V\bfx)$.
Using inequality (2.1) in \cite{KobayashiTsuchiya3}, we have
\begin{gather*}
  \frac{\prod_{i=1}^2(1-\bfs_i)^{r}}{t_1^{2r}t_2^{2r}}
  \sum_{|\gamma|=r} (\partial_\bfx^\gamma u)^2
 \le \sum_{|\gamma|=r} (\partial_\bfy^\gamma v)^2 \le
\frac{\prod_{i=1}^2(1+\bfs_i)^{r}}{t_1^{2r}t_2^{2r}}
  \sum_{|\gamma|=r} (\partial_\bfx^\gamma u)^r,
\end{gather*}
where $\bfy := V\bfx$.  The fact that
$\mathrm{det}A = t_1t_2$ and inequalities (2.2), (2.3) in
\cite{KobayashiTsuchiya3} then give
\begin{align*}
   |v|_{m,p,K} & \le 3^{m\tau(p)} 
 \frac{\prod_{i=1}^2(1+\bfs_i)^{m/2}}
       {(t_1 t_2)^{m-1}}  |u|_{m,p,\Kabc}, \\
  |v|_{k+1,p,K} & \ge 3^{-(k+1)\tau(p)}   
\frac{\prod_{i=1}^2(1-\bfs_i)^{(k+1)/2}}{(t_1t_2)^{k}}
  |u|_{k+1,p,\Kabc},   \\
   \frac{|v|_{m,p,K}}{|v|_{k+1,p,K}} & \le 3^{(k+1+m)\tau(p)}
   \frac{(t_1t_2)^{k+1-m}\prod_{i=1}^2(1+\bfs_i)^{m/2}|u|_{m,p,\Kabc}}
     {\prod_{i=1}^2(1-\bfs_i)^{(k+1)/2}|u|_{k+1,p,\Kabc}} \\
 & = 3^{(k+1+m)\tau(p)} \frac{\prod_{i=1}^2(1+\bfs_i)^{(k+1)/2}}
  {\prod_{i=1}^2(1-\bfs_i)^{m/2}}
  \frac{|u|_{m,p,\Kabc}}{|u|_{k+1,p,\Kabc}} \\
  & \le C  \frac{\left(\max\{\alpha,\gamma\}\right)^{k+1-m}}
  {\prod_{i=1}^2(1-\bfs_i)^{m/2}},
\end{align*}
where
\begin{align*}
  \tau(p) := \begin{cases}
    \frac{1}{p} - \frac{1}{2}, & 1 \le p \le 2 \\
    \frac{1}{2} - \frac{1}{p} , & 2 \le p \le \infty
      \end{cases}, \qquad
  C := 3^{(k+1+m)\tau(p)} 2^{k+1} C_{k,m,p}.
\end{align*}
(See the inequalities in \cite[p.496]{KobayashiTsuchiya3}.)
Note that the constant $C$ depends only on $k$, $m$, and $p$.
Hence, we have derived the following theorem:
\begin{theorem}\label{thmB}
Let $K$ be an arbitrary tetrahedron with vertices given by
\eqref{eq:position} with the parameters in \eqref{eq:tetra-param}.
Let $k$, $m$ be integers with $k \ge 1$ and $0 \le m \le k$.
Let $p$ be taken as \eqref{p-cond}.  Then, we have
\[
 B_p^{m,k}(K) := \sup_{v \in \T_p^k}\frac{|v|_{m,p,K}}{|v|_{k+1,p,K}}
  \le C  \frac{\left(\max\{\alpha,\gamma\}\right)^{k+1-m}}
  {\prod_{i=1}^2(1-\bfs_i)^{m/2}},
\]
where $C = C(k,m,p)$ is a constant independent of $K$.
\end{theorem}

\section{A geometric interpretation and the proof of the main theorem}
\label{sec:geometric}
In this section, we prove the main theorem (Theorem~\ref{thm:main})
from Theorem~\ref{thmB}. 
To this end, we consider the geometric meaning of the quantity
$\prod_{i=1}^2(1-\bfs_i)^{-1/2}$ that appeared in Theorem~\ref{thmB}.
Recall that $K$ is a tetrahedron with vertices given by
\eqref{eq:position} and the parameters in \eqref{eq:tetra-param}.
Recall also that $B$ is the base of $K$ with vertices
$\bfx_1$, $\bfx_2$, $\bfx_3$.  The circumradius $R_B$
can then be written as
\[
   R_B = \frac{\sqrt{\alpha^2 - 2\alpha\beta s_1 + \beta^2}}{2t_1}.
\]
Because
\begin{align*}
   \alpha^2 - 2\alpha\beta s_1 + \beta^2 & =
  \frac{\alpha^2(2 - s_1^2)}{4} + \alpha
   \left(\frac{\alpha}{2} - \beta s_1\right) 
   + \left(\frac{\alpha s_1}{2} - \beta\right)^2 \\
   & \ge \frac{\alpha^2(2 - s_1^2)}{4} \ge \frac{\alpha^2}{4},
\end{align*}
we have
\begin{align}
   R_B \ge \frac{\alpha}{4 t_1} = \frac{\alpha}{4 \sqrt{1 - \bfs_1^2}} 
    \ge \frac{h_B}{4\sqrt{2} \sqrt{1 - \bfs_1}},
   \label{RB-est}
\end{align}
where we have used the definition $h_B = \alpha$.

Recall that $R_P$ was defined in Section~\ref{section:Circumradius}.
We will show that there exists a constant $C$ independent of $K$
such that
\begin{align}
    R_P \ge C\frac{\max\{\alpha,\gamma\}}{\sqrt{1 - \bfs_2}}.
   \label{RP-est}
\end{align}

Let $\bfx_3 = (\eta,\xi,0)^\top$, that is, $\xi = \beta t_1$ and
$\eta = \beta s_1$ or $\eta = \alpha - \beta s_1$.  From the
assumption, we have $0 < \eta < \alpha$ and $\xi > 0$.  Note that 
\begin{align*}
   \delta_\theta(\bfx_1) = (0,0,0)^\top, \quad
  & \delta_\theta(\bfx_2) = (\alpha\cos\theta,0,0)^\top, \\
   \delta_\theta(\bfx_3) = (\eta\cos\theta-\xi\sin\theta,0,0)^\top, 
  \quad &
   \delta_\theta(\bfx_4) = \gamma
  (s_{21}\cos\theta-s_{22}\sin\theta,0,t_2)^\top.
\end{align*}

Let $\ux < \ox$ be the $x$-coordinates of the end points of the base of
$\delta_\theta(K)$.  The assumptions given in \eqref{eq:tetra-param}
yield
\begin{align*}
   \ux \le 0, \quad \ox = \alpha\cos\theta, \quad 
   0 \le \theta \le \frac{\pi}{2}, \qquad
   \ux =0, \quad \ox \ge \alpha\cos\theta, \quad 
   -\frac{\pi}{2} \le \theta \le 0. 
\end{align*}
Defining
$w = w(\theta) := s_{21}\cos\theta - s_{22}\sin\theta$, $R_\theta$ is
written as
\[
   R_\theta = \frac{1}{2\gamma t_2}
   \left((\ox - \gamma w)^2 + \gamma^2t_2^2\right)^{1/2}
   \left((\ux - \gamma w)^2 + \gamma^2t_2^2\right)^{1/2}.
\]
Take an arbitrary
$\theta \in \left[-\frac{\pi}{3},\frac{\pi}{3}\right]$.
Suppose that the inequality
\begin{align}
  \gamma w \le \frac{\ux + \ox}{2}
  \label{key-inequality1}
\end{align}
holds and there exists a constant $C_1$ independent of
$\theta$ such that
\begin{align}
 |\ux - \gamma w| \ge C_1 \gamma \bfs_2.  \label{key-inequality2}
\end{align}
We then have
\begin{align*}
   (\ux - \gamma w)^2 + \gamma^2 t_2^2 & \ge C_1^2
   (1 - t_2^2) + \gamma^2 t_2^2 \ge
   \min\{1, C_1^2\} \gamma^2, \\
   (\ox - \gamma w)^2 + \gamma^2 t_2^2 & \ge \max\left\{
      (\ux - \gamma w)^2 + \gamma^2 t_2^2, \;
      \frac{\alpha^2\cos^2\theta}{4} + \gamma^2 t_2^2 \right\} \\
  & \ge \max\left\{\min\{1,C_1^2\} \gamma^2,\, C_2^2\alpha^2\right\} \\
  & \ge \min\{C_1^2,C_2^2\}\max\{\alpha^2,\gamma^2\},
\end{align*}
where $C_2 = C_2(\theta) := \frac{\cos\theta}{2}$.
Here, we have used the fact that $\ox/2 \ge (\ox+\ux)/2 \ge \gamma w$ and
$\ox - \gamma w \ge \ox/2 \ge (\alpha\cos\theta)/2$.  Hence,
setting
\[
   \frac{\min\{C_1,C_2\}\min\{1,C_1\}}{2\sqrt{2}} 
  \ge \frac{\min\left\{C_1,\frac{1}{4}\right\}
   \min\{1,C_1\}}{2\sqrt{2}} =: C_3,
\]
we obtain
\begin{align*}
  R_P \ge R_\theta \ge \frac{\sqrt{2}C_3}{t_2}\max\{\alpha,\gamma\}
   = \frac{\sqrt{2}C_3 \max\{\alpha,\gamma\}}{\sqrt{1-\bfs_2^2}}
   \ge C_3 \frac{\max\{\alpha,\gamma\}}{\sqrt{1-\bfs_2}},
\end{align*}
and the key inequality \eqref{RP-est} is shown.

Fix $\varphi$ such that
\[
   0 < \varphi < \frac{\pi}{6}, \qquad 
   \sin2\varphi\tan2\varphi \le \frac{1}{6}.
\]
In the following, we will show that, according to 
$\bfx_4 = (\gamma s_{21}, \gamma s_{22}, \gamma t_2)^\top$, we can
take an appropriate
$\theta \in \left[-\frac{\pi}{3},\frac{\pi}{3}\right]$ such that
conditions \eqref{key-inequality1} and \eqref{key-inequality2}
hold with $C_1 = \sin\varphi$.  

\vspace{0.3cm}
\noindent
\textbf{Case 1.} Suppose that $|s_{22}|\tan\varphi \le |s_{21}|$. \\
In this case, we set $\theta = 0$ and have $\ux = 0$, $\ox = \alpha$,
and $\gamma w = \gamma s_{21} \le \alpha/2 = (\ux + \ox)/2$
because of \eqref{eq:tetra-param}.  Hence, \eqref{key-inequality1} holds.
For \eqref{key-inequality2}, we note that
\begin{align*}
   |\ux - \gamma w| = \gamma|s_{21}| & = \gamma
   \left(s_{21}^2\sin^2\varphi + s_{21}^2\cos^2\varphi\right)^{1/2} \\
  & \ge \gamma
   \left(s_{21}^2\sin^2\varphi + s_{22}^2\sin^2\varphi\right)^{1/2}
   = \gamma \bfs_2 \sin\varphi,
\end{align*}
and so \eqref{key-inequality2} holds with $C_1 := \sin\varphi$.

\vspace{0.3cm}
\noindent
\textbf{Case 2:} Suppose that $|s_{22}|\tan\varphi > |s_{21}|$
and $3\gamma s_{22} \tan2\varphi \le \alpha$. \\
In this case, set $\theta = -2\varphi$. We then have $\ux = 0$
and $\ox \ge \alpha\cos\theta = \alpha \cos2\varphi$.  If
$s_{22} > 0$, then
\begin{align*}
 \gamma w & = \gamma (s_{21}\cos2\varphi + s_{22}\sin2\varphi)
  \le \gamma (s_{22}\cos2\varphi\tan\varphi + s_{22}\sin2\varphi) \\
  & = \frac{\gamma s_{22}}{2}
   (3 \sin2\varphi - 2\sin^2\varphi\tan\varphi) \\
  & \le \frac{2\gamma s_{22}}{2}\sin2\varphi
   \le \frac{\alpha \cos2\varphi}{2} \le \frac{\ux + \ox}{2}.
\end{align*}
   If $s_{22} \le 0$, then
\begin{align*}
 \gamma w & = \gamma (s_{21}\cos2\varphi + s_{22}\sin2\varphi)
  \le \gamma (- s_{22}\cos2\varphi\tan\varphi + s_{22}\sin2\varphi) \\
  & = \gamma s_{22}\tan\varphi \le 0 \le \frac{\ux + \ox}{2}.
\end{align*}
Thus, in either case, \eqref{key-inequality1} holds.
For \eqref{key-inequality2}, we note that
\begin{align*}
   |\ux - \gamma w| & = \gamma
  |s_{21}\cos2\varphi + s_{22}\sin2\varphi| \ge \gamma \left(
  |s_{22}|\sin2\varphi - |s_{21}|\cos2\varphi\right) \\
  & \ge \gamma \left(
  |s_{22}|\sin2\varphi - |s_{22}|\cos2\varphi\tan\varphi\right)
  = \gamma|s_{22}|\tan\varphi \\
  & = \gamma \left(
    s_{22}^2 \sin^2\varphi\tan^2\varphi + s_{22}^2\sin^2\varphi
   \right)^{1/2} \\
  & = \gamma \left(s_{21}^2 \sin^2\varphi + s_{22}^2\sin^2\varphi
   \right)^{1/2} = \gamma\bfs_2\sin\varphi,
  \end{align*}
and so \eqref{key-inequality2} holds with $C_1 := \sin\varphi$.

\vspace{0.3cm}
\noindent
\textbf{Case 3:} Suppose that $|s_{22}|\tan\varphi > |s_{21}|$
and $3\gamma s_{22} \tan2\varphi > \alpha$. \\
In this case, set $\theta = 2\varphi$.  We then have
\[
   \ux = \min\{\eta\cos2\varphi - \xi\sin2\varphi,0\}, \qquad
    \ox = \alpha \cos2\varphi.
\]
If $\eta\cos2\varphi - \xi\sin2\varphi \le 0$, then
\begin{align*}
   \ux & = \eta\cos2\varphi - \xi\sin2\varphi
   = \alpha \cos2\varphi + (\eta - \alpha) \cos2\varphi
      - \xi\sin2\varphi \\
  & \ge \alpha \cos2\varphi -
   \left((\eta - \alpha)^2 + \xi^2\right)^{1/2} 
   \ge \alpha \cos2\varphi - \alpha = -2\alpha \sin^2\varphi,
\end{align*}
because the lengths of all edges of the base are less than $\alpha$.
Even if $\eta\cos2\varphi - \xi\sin2\varphi \ge 0 = \ux$, the above
inequality obviously holds.  Because 
\begin{align*}
  \gamma w & = \gamma (s_{21}\cos2\varphi - s_{22}\sin2\varphi) \\
  & \le \gamma (s_{22}\cos2\varphi\tan\varphi - s_{22}\sin2\varphi)
  = -\gamma s_{22}\tan\varphi,
\end{align*}
we have
\begin{align*}
  \ux - \gamma w & \ge -2\alpha\sin^2\varphi + \gamma s_{22}\tan\varphi \\
  & \ge -2\alpha\sin^2\varphi + \frac{\alpha\tan\varphi}{6\tan2\varphi}
 + \frac{\gamma}{2} s_{22}\tan\varphi \\
  & = \alpha\left(\frac{1}{6}- \sin2\varphi\tan2\varphi\right)
   \frac{\tan\varphi}{\tan2\varphi}
   + \frac{\gamma}{2} s_{22}\tan\varphi
    \ge \frac{\gamma}{2} s_{22}\tan\varphi > 0.
\end{align*}
Therefore,
\[
   \gamma w < \ux < \frac{\ux + \ox}{2},
\]
and \eqref{key-inequality1} holds.
For \eqref{key-inequality2}, we note that
\begin{align*}
   |\ux - \gamma w| & \ge \frac{\gamma}{2} s_{22}\tan\varphi \\
  & = \gamma \left(
    s_{22}^2 \sin^2\varphi\tan^2\varphi + s_{22}^2\sin^2\varphi
   \right)^{1/2} \\
  & \ge \gamma \left(s_{21}^2 \sin^2\varphi + s_{22}^2\sin^2\varphi
   \right)^{1/2} = \gamma\bfs_2\sin\varphi,
  \end{align*}
and so \eqref{key-inequality2} holds with $C_1 := \sin\varphi$.

\vspace{0.3cm}
Using inequalities \eqref{RB-est} and \eqref{RP-est}, we have shown the
following lemma.

\begin{lemma}\label{lemma:geometric}
Let $K$ be a tetrahedron with vertices given by \eqref{eq:position}
and the parameters in \eqref{eq:tetra-param}.
We then have
\[
  \prod_{i=1}^2(1-\bfs_i)^{-1/2} \le C
   \frac{R_BR_P}{h_B\max\{\alpha,\gamma\}},
\]
where $C$ is a constant independent of $K$.
\end{lemma}

Take a facet $B$ so that $R_BR_P/h_B = R_K$.
Combining Theorem~\ref{thmB} and Lemma~\ref{lemma:geometric} with the 
projected circumradius $R_K$, we have
\begin{align*}
 B_p^{m,k}(K) := \sup_{v \in \T_p^k}\frac{|v|_{m,p,K}}{|v|_{k+1,p,K}}
  & \le C  \frac{\left(\max\{\alpha,\gamma\}\right)^{k+1-m}}
  {\prod_{i=1}^2(1-\bfs_i)^{m/2}} \\
  & \le C \left(\frac{R_BR_P}{h_B}\right)^m
    \left(\max\{\alpha,\gamma\}\right)^{k+1-2m} \\
  & \le C R_K^m h_K^{k+1-2m},
\end{align*}
where the constant $C = C(k,m,p)$ is independent of $K$
with vertices \eqref{eq:position}.  Note that Sobolev (semi-)norms 
are affected by rotation up to a constant.  Therefore, we have
proved the main theorem (Theorem~\ref{thm:main}).

\section{Concluding remarks}
In this final section, we present some remarks on the newly obtained
error estimation.

(1)  Let $\Omega \subset \R^3$ be a bounded polyhedral domain.
Suppose that we compute a numerical solution of the Poisson equation
\begin{align}
    -\Delta u = f \text{ in } \Omega, \quad
    u = 0 \text{ on } \partial\Omega
   \label{model-eq}
\end{align}
by the piecewise $k$th-order finite element method with conforming
simplicial elements.  To this end, we construct a tetrahedralization
$\T_h$ of $\Omega$ and consider the piecewise $\PP_k$ continuous
function space $S_h \subset H_0^1(\Omega)$.
Then, the weak form of \eqref{model-eq} is
\begin{align*}
    \int_\Omega \nabla u \cdot \nabla v \dd \bfx  = 
    \int_\Omega f v \dd \bfx, \quad \forall
      v \in H_0^1(\Omega),
\end{align*}
and the finite element solution is defined as the unique solution
$u_h \in S_h$ of
\begin{align*}
    \int_\Omega \nabla u_h \cdot \nabla v_h \dd \bfx  = 
    \int_\Omega f v_h \dd \bfx, \quad \forall v_h \in S_h.
\end{align*}
C\'ea's Lemma implies that the error $|u - u_h|_{1,2,\Omega}$ is
estimated as
\begin{align}
  |u - u_h|_{1,2,\Omega} \le \left(
   \sum_{K \in \T_h} |u - \I_K^k u|_{1,2,K}^2 \right)^{1/2}.
  \label{Cea_lemma}
\end{align}
Combining \eqref{Cea_lemma} and Theorem~\ref{thm:main} with $p=2$,
$k \ge 2$, $m=1$, we have
\begin{align*}
  |u - u_h|_{1,2,\Omega} & \le C\left(\sum_{K \in \T_h}
   |u - \I_K^ku|_{k+1,2,K}^2\right)^{1/2} \\
  & \le C\left(\sum_{K \in \T_h}
   (R_K h_K^{k-1})^2 |u|_{k+1,2,K}^2\right)^{1/2} \\
  & \le C \max_{K \in \T_h}(R_Kh_K^{k-1}) |u|_{k+1,2,\Omega}. 
\end{align*}
Therefore, if $\max_{K \in \T_h}(R_Kh_K^{k-1}) \to 0$ as $h \to 0$ and
$u \in W^{k+1,2}(\Omega)$, the finite element solution $u_h$ converges
to the exact solution $u$ even if there exist many skinny elements
violating the maximum angle condition in $\T_h$.

It is known that a triangulation of 3-dimensional domain $\Omega$
may have many slivers.  Consider the typical sliver $K$ mentioned in
Section~\ref{intro}, whose vertices are
$(\pm h,0,0)^\top$ and $(0,\pm h,h^\alpha)^\top$.  We can see that
$R_K = \mathcal{O}(h_K^{2-\alpha})$ and
$R_Kh_K^{k-1} = \mathcal{O}(h_K^{k+1-\alpha})$.  Thus, if we use
a $k$th-order conforming Lagrange finite element method with
triangulation that contains many slivers like $K$, the theoretical
convergence rate may be $\mathcal{O}(h^{k+1-\alpha})$.  This is
worse than the expected rate $\mathcal{O}(h^{k})$, but we can still
expect convergence if $k+1-\alpha > 0$.   Therefore, ``bad''
triangulation with many slivers can be remedied by using
higher-order Lagrange elements.

(2) In Theorem~\ref{thmB}, we obtained an error estimation for
Lagrange interpolation on tetrahedrons in terms of the singular values
${\prod_{i=1}^2(1-\bfs_i)^{-1/2}}$ of the linear transformation.
In Lemma~\ref{lemma:geometric}, we showed that the projected circumradius
is a geometric interpretation of the singular values.  The authors,
however, are not completely sure whether the projected circumradius is
the \textit{best} interpretation. 
Further research on the geometry of tetrahedrons is required 
to ascertain the geometric properties of the singular values 
${\prod_{i=1}^2(1-\bfs_i)^{-1/2}}$ and the projected circumradius.

(3) The finite element error analysis on triangulations with
anisotropic elements have been done by many authors in many ways.
See, for example, \cite{ApelDobrowolski}, \cite{ApelLube}, \cite{Apel}.
 \cite{FormaggiaPerotto}, and \cite{Cao}.
The present authors believe that the approach and results given in
this paper are fundamentally different from those papers and give
a new insight to the finite element error analysis.

\section{Appendix:  A counter example of the Squeezing Theorem with
 $p=2$ and $k=m=1$}
In this Appendix, we present a counter example which shows that the
Squeezing Theorem (Theorem~\ref{squeezingtheorem}) and thus the main
theorem (Theorem~\ref{thm:main}) do not hold for the case $p=2$ and
$k=m=1$.  Our counter example is simpler than Shenk's counter example
given in \cite{Shenk}.  Let $\hK_h$ be the tetrahedron with vertices
$(0,0,0)^\top$, $(1,0,0)^\top$, $(0,1,0)^\top$, and $(0,0,h)^\top$.  For
$k > 0$, define
\[
  u(x,y,z) := z \log(1 + 2k(x+y)). 
\]
Clearly, $u$ vanishes at the vertices of $\hK_h$.
By straightforward computation, we find
{\allowdisplaybreaks
\begin{align*}
  |u|_{1,2,\hK_h}^2 & = \int_{\hK_h}
  \left(|u_x|^2 + |u_y|^2 + |u_z|^2\right) \dd \bfx  
  \ge \int_{\hK_h} |u_z|^2 \dd \bfx \\
  & = \int_0^1 \int_0^{1-x} \int_0^{h(1-x-y)}
   \left\{\log (1 + 2k(x+y))\right\}^2 \dd z \dd y \dd x \\
  & = h \int_0^1 t(1-t) \left(\log (1 + 2kt)\right)^2 \dd t \\
  & \ge h \int_0^{1/2} t(1-t) \left(\log (1 + 2kt)\right)^2 \dd t 
    \ge \frac{h}{2} \int_0^{1/2} t \left(\log (1 + 2kt)\right)^2 \dd t \\
    & \ge \frac{h}{4} \int_0^{1/2} t \log (1 + 2kt)
       \left(\log (1 + 2kt) + \frac{2kt}{1+2kt}\right) \dd t \\
  & = \frac{h}{8}\left[t^2 \left(\log(1+2kt)\right)^2\right]_0^{1/2}
    = \frac{h}{32}(\log(1+k))^2,
\end{align*}
}
in which we have used the facts 
\[
   \int_0^1 \int_0^{1-x} f(x+y) \dd y \dd x = \int_0^1 tf(t) \dd t
  \quad \text{ and } \quad \log(1 + x) \ge \frac{x}{1+x}, \quad x \ge 0.
\]
We also find
{\allowdisplaybreaks
\begin{align*}
 |u|_{2,2,\hK_h}^2 & = \int_0^1 \int_0^{1-x} \int_0^{h(1-x-y)}
   \left(\frac{64k^4z^2}{(1+2k(x+y))^4}
     + \frac{16k^2}{(1+2k(x+y))^2} \right) \dd z \dd y \dd x \\
  & \le \int_0^1 \int_0^{1-x} \int_0^{h}
   \left(\frac{64k^4 h^2}{(1+2k(x+y))^4}
     + \frac{16k^2}{(1+2k(x+y))^2} \right) \dd z \dd y \dd x \\
  & \le 8 k h \int_0^1 \int_0^{1-x} \left(\frac{8k^3 h^2}{(1+2k(x+y))^4}
     + \frac{2k}{(1+2k(x+y))^2} \right) \dd y \dd x \\
  & = 8 k h \int_0^1 \left(\frac{8k^3 h^2 t}{(1+2kt)^4}
     + \frac{2kt}{(1+2kt)^2} \right) \dd t \\
  & \le 8 k h \int_0^1 \left(\frac{4k^2 h^2}{(1+2kt)^3}
     + \frac{1}{1+2kt} \right) \dd t \\
  & \le 8 k h \int_0^\infty \frac{4k^2 h^2}{(1+2kt)^3} \dd t
     + 8 k h \int_0^1 \frac{1}{1+kt} \dd t
  = 8h (k^2h^2 + \log(1+k)),
\end{align*}
}
in which we have used the fact $x/(1+x) \le 1$, $x \ge 0$. 
This shows that $u \in H^2(\hK_h)$ and $u \in \T_2^1(\hK_h)$.
Combining these two inequalities, we have
\begin{align*}
  \frac{|u|_{1,2,\hK_h}^2}{|u|_{2,2,\hK_h}^2}
  \ge \frac{(\log(1+k))^2}{256(k^2h^2 + \log(1+k))}.
\end{align*}
Therefore, setting $k:=1/h$, we find
\begin{align*}
  \frac{|u|_{1,2,\hK_h}^2}{|u|_{2,2,\hK_h}^2}
  \ge \frac{(\log(1+1/h))^2}{256(1 + \log(1+1/h))} \to \infty
 \quad \text{ as } \quad h \to 0.
\end{align*}
This means that 
\begin{align*}
 B_2^{1,1}(\hK_h) := \sup_{v \in \T_2^1(\hK_h)}
  \frac{|v|_{1,2,\hK_h}}{|v|_{2,2,\hK_h}}
   \to \infty  \quad \text{ as } \quad h \to 0,
\end{align*}
and Theorem~\ref{squeezingtheorem} does not hold for the case
$p=2$ and $k=m=1$.

\vspace{0.5cm}
\noindent
\textbf{Acknowledgments}
The authors were supported by JSPS KAKENHI Grant Numbers
 JP26400201, JP16H03950, and JP17K18738.
The authors thank the anonymous referees for their
valuable comments.

\section*{References}


\end{document}